# Self-dual Modules of Semisimple Hopf Algebras

Yevgenia Kashina    Yorck Sommerhäuser

Yongchang Zhu

**Abstract**

We prove that, over an algebraically closed field of characteristic zero, a semisimple Hopf algebra that has a nontrivial self-dual simple module must have even dimension. This generalizes a classical result of W. Burnside. As an application, we show under the same assumptions that a semisimple Hopf algebra that has a simple module of even dimension must itself have even dimension.

**1**  Suppose that $H$ is a finite-dimensional Hopf algebra that is defined over the field $K$. We denote its comultiplication by $\Delta$, its counit by $\varepsilon$, and its antipode by $S$. For the comultiplication, we use the sigma notation of R. G. Heyneman and M. E. Sweedler in the following variant:

$$\Delta(h) = h_{(1)} \otimes h_{(2)}$$

We view the dual space $H^*$ as a Hopf algebra whose unit is the counit of $H$, whose counit is the evaluation at 1, whose antipode is the transpose of the antipode of $H$, and whose multiplication and comultiplication are determined by the formulas

$$(\varphi\varphi')(h) = \varphi(h_{(1)})\varphi'(h_{(2)}) \qquad \varphi_{(1)}(h)\varphi_{(2)}(h') = \varphi(hh')$$

for $h, h' \in H$ and $\varphi, \varphi' \in H^*$.

With $H$, we can associate its Drinfel'd double $D(H)$ (cf. [18], § 10.3, p. 187). This is a Hopf algebra whose underlying vector space is $D(H) = H^* \otimes H$. As a coalgebra, it is the tensor product of $H^{*\,\mathrm{cop}}$ and $H$, i.e., we have

$$\Delta(\varphi \otimes h) = (\varphi_{(2)} \otimes h_{(1)}) \otimes (\varphi_{(1)} \otimes h_{(2)})$$

as well as $\varepsilon(\varphi \otimes h) = \varphi(1)\varepsilon(h)$. Its multiplication is given by the formula

$$(\varphi \otimes h)(\varphi' \otimes h') = \varphi'_{(1)}(S^{-1}(h_{(3)}))\varphi'_{(3)}(h_{(1)})\varphi\varphi'_{(2)} \otimes h_{(2)}h'$$

The unit element is $\varepsilon \otimes 1$ and the antipode is $S(\varphi \otimes h) = (\varepsilon \otimes S(h))(S^{*-1}(\varphi) \otimes 1)$.

**2** As the underlying vector space of $D(H)$ is $H^* \otimes H$, there is a canonical linear form on $D(H)$, namely the evaluation form:

$$e : D(H) \to K, \ \varphi \otimes h \mapsto \varphi(h)$$

This form is an invertible element of $D(H)^*$; its inverse is given by the formula $e^{-1}(\varphi \otimes h) = \varphi(S^{-1}(h))$. This holds since

$$e^{-1}(\varphi_{(2)} \otimes h_{(1)})e(\varphi_{(1)} \otimes h_{(2)}) = \varphi_{(2)}(S^{-1}(h_{(1)}))\varphi_{(1)}(h_{(2)}) = \varphi(1)\varepsilon(h)$$

A similar calculation shows that $e^{-1}$ is also a right inverse of $e$.

The evaluation form was considered by T. Kerler, who proved the following property (cf. [12], Prop. 7, p. 366):

**Proposition 1** The evaluation form is a symmetric Frobenius homomorphism.

**Proof.** We give a different proof. By the definition of a Frobenius algebra (cf. [10], Kap. 13, Def. 13.5.4, p. 306), we have to show that the bilinear form associated with $e$ is symmetric and nondegenerate. Since we have

$$\begin{aligned} e((\varphi \otimes h)(\varphi' \otimes h')) &= \varphi'_{(1)}(S^{-1}(h_{(3)}))\varphi'_{(3)}(h_{(1)}) \ e(\varphi\varphi'_{(2)} \otimes h_{(2)}h') \\ &= \varphi'_{(1)}(S^{-1}(h_{(4)}))\varphi'_{(3)}(h_{(1)})\varphi(h_{(2)}h'_{(1)})\varphi'_{(2)}(h_{(3)}h'_{(2)}) \\ &= \varphi'(S^{-1}(h_{(4)})h_{(3)}h'_{(2)}h_{(1)})\varphi(h_{(2)}h'_{(1)}) = \varphi'(h'_{(2)}h_{(1)})\varphi(h_{(2)}h'_{(1)}) \end{aligned}$$

we see that this bilinear form is symmetric. To see that it is also nondegenerate, consider the right multiplication $R_e$ by $e$ in $D(H)^*$. By dualizing this map, we get the following endomorphism of $D(H)$:

$$R_e^* : D(H) \to D(H), \ \varphi \otimes h \mapsto \varphi_{(1)}(h_{(2)})\varphi_{(2)} \otimes h_{(1)}$$

The inverse of this endomorphism is obviously obtained by dualizing the right multiplication by $e^{-1}$:

$$R_{e^{-1}}^* : D(H) \to D(H), \ \varphi \otimes h \mapsto \varphi_{(1)}(S^{-1}(h_{(2)}))\varphi_{(2)} \otimes h_{(1)}$$

Since from the above we have that

$$\begin{aligned} &e(R_{e^{-1}}^*(\varphi \otimes h)R_{e^{-1}}^*(\varphi' \otimes h')) \\ &= \varphi_{(1)}(S^{-1}(h_{(2)}))\varphi'_{(1)}(S^{-1}(h'_{(2)}))e((\varphi_{(2)} \otimes h_{(1)})(\varphi'_{(2)} \otimes h'_{(1)})) \\ &= \varphi_{(1)}(S^{-1}(h_{(3)}))\varphi'_{(1)}(S^{-1}(h'_{(3)}))\varphi'_{(2)}(h'_{(2)}h_{(1)})\varphi_{(2)}(h_{(2)}h'_{(1)}) = \varphi(h')\varphi'(h) \end{aligned}$$

the bilinear form under consideration is isometric to the bilinear form

$$D(H) \times D(H) \to K, \ (\varphi \otimes h, \varphi' \otimes h') \mapsto \varphi(h')\varphi'(h)$$

which is obviously nondegenerate. $\square$

The powers of $e$ are given by the following formula:

$$e^m(\varphi \otimes h) = e(\varphi_{(m)} \otimes h_{(1)})e(\varphi_{(m-1)} \otimes h_{(2)}) \cdot \ldots \cdot e(\varphi_{(1)} \otimes h_{(m)}) = \varphi(h_{(m)} \cdot \ldots \cdot h_{(1)})$$

This shows that the order of $e$ is related to the exponent of $H$:

**Proposition 2** Suppose that $H$ is semisimple and that the base field $K$ has characteristic zero. Then the order of $e$ is equal to the exponent of $H$. In particular, the order of $e$ divides $(\dim(H))^3$.

**Proof.** In this situation, we know from [15], Thm. 3.3, p. 276, and [16], Thm. 3, p. 194 that $H$ is also cosemisimple and that the antipode of $H$ is an involution. It therefore follows from the definition of the exponent (cf. [5], Def. 2.1, p. 132) that the order of $e$ is the exponent of $H^{\text{op}}$, which coincides with the exponent of $H$ by [5], Cor. 2.6, p. 134. The divisibility property is proved in [5], Thm. 4.3, p. 136. $\square$

**3** Let us consider now the case that $H$ is semisimple and that the base field $K$ is algebraically closed of characteristic zero. Note that a semisimple Hopf algebra is necessarily finite-dimensional (cf. [22], Cor. 2.7, p. 330, or [23], Chap. V, Ex. 4, p. 108). By Maschke's theorem (cf. [18], Thm. 2.2.1, p. 20), there is a unique two-sided integral $\Lambda$ that satisfies $\varepsilon(\Lambda) = 1$. Suppose that $V$ is a simple $H$-module with character $\chi$. We say that $V$ is self-dual if $V \cong V^*$. This is equivalent to the requirement that there is a nondegenerate invariant bilinear form on $V$, i.e., a nondegenerate bilinear form

$$\langle \cdot, \cdot \rangle : V \times V \to K$$

that satisfies

$$\langle h_{(1)}.v, h_{(2)}.v' \rangle = \varepsilon(h) \langle v, v' \rangle$$

for all $h \in H$ and all $v, v' \in V$. Following [17], we define the Frobenius-Schur indicator, also briefly called the Schur indicator, $\nu_2(\chi)$ of the irreducible character $\chi$ corresponding to the simple module $V$:

$$\nu_2(\chi) := \chi(\Lambda_{(1)} \Lambda_{(2)})$$

The Frobenius-Schur theorem for Hopf algebras (cf. [17], Thm. 3.1, p. 349) then asserts, among other things, the following:

**Theorem** The Schur indicator $\nu_2(\chi)$ can only take the values $1$, $-1$, and $0$:

1. We have $\nu_2(\chi) = 1$ if and only if $V$ admits a symmetric nondegenerate invariant bilinear form.

2. We have $\nu_2(\chi) = -1$ if and only if $V$ admits a skew-symmetric nondegenerate invariant bilinear form.

3. We have $\nu_2(\chi) = 0$ if and only if $V$ is not self-dual.

**4** Using these preparations, we can prove the main theorem. It generalizes a classical result of W. Burnside in the theory of finite groups (cf. [3], Par. 2, p. 167; [4], § 222, Thm. II, p. 294). We note that this theorem was known in the case of cocentral abelian extensions (cf. [11], Cor. 3.2, p. 5).

**Theorem** Suppose that $H$ is a semisimple Hopf algebra over an algebraically closed field of characteristic zero. If $H$ has a nontrivial self-dual simple module, then the dimension of $H$ is even.

**Proof.** Suppose that $V$ is an $H$-module with character $\chi$ and that $W$ is an $H^*$-module with character $\eta$. As an algebra, the dual $D(H)^*$ of the Drinfel'd double is isomorphic to $H^{\mathrm{op}} \otimes H^*$. We can therefore turn $V \otimes W$ into a $D(H)^*$-module by defining
$$(h \otimes \varphi).(v \otimes w) = S(h).v \otimes \varphi.w$$
If we identify $H^{**}$ and $H$, we can consider $\eta$ as an element of $H$. Denoting the character of $V^*$ by $\bar{\chi}$, the trace of the action of $e$ on $V \otimes W$ is then given by the formula
$$(\bar{\chi} \otimes \eta)(e) = e(\bar{\chi} \otimes \eta) = \bar{\chi}(\eta)$$
Similarly, the trace of $e^2$ is given by the formula
$$(\bar{\chi} \otimes \eta)(e^2) = e^2(\bar{\chi} \otimes \eta) = \bar{\chi}(\eta_{(2)}\eta_{(1)})$$
We now assume that $V$ is simple, nontrivial, and self-dual and that $W = H^*$ is the regular representation. We then know that, if $\Lambda$ is an integral that satisfies $\varepsilon(\Lambda) = 1$, the character of the regular representation is given by
$$\eta(\varphi) = (\dim H)\varphi(\Lambda)$$
i.e., up to the identification of $H^{**}$ and $H$, we have $\eta = (\dim H)\Lambda$. Since $V$ is nontrivial, $\chi$ vanishes on the integral, and since the self-duality of $V$ implies that $\bar{\chi} = \chi$, we get from the above and the Frobenius-Schur theorem that
$$(\bar{\chi} \otimes \eta)(e) = 0 \qquad (\bar{\chi} \otimes \eta)(e^2) = \pm \dim(H)$$

Now suppose that $n$ is the exponent of $H$ and that $\zeta$ is a primitive $n$-th root of unity. Since $e$ has order $n$ by Proposition 2.2, $V \otimes W$ is the direct sum of the eigenspaces corresponding to the powers of $\zeta$, whose dimensions we denote by
$$a_k := \dim\{z \in V \otimes W \mid e.z = \zeta^k z\}$$
If we introduce the polynomial
$$p(x) := \sum_{k=0}^{n-1} a_k x^k \in \mathbb{Z}[x]$$

we see that $p(\zeta) = (\bar{\chi} \otimes \eta)(e) = 0$. Therefore, if $q_n$ denotes the $n$-th cyclotomic polynomial, we see that $q_n$ divides $p$. On the other hand, $e^2$ acts on the eigenspace of $e$ corresponding to the eigenvalue $\zeta^i$ by multiplication with $\zeta^{2i}$. Therefore, we get

$$p(\zeta^2) = (\bar{\chi} \otimes \eta)(e^2) = \pm \dim(H) \neq 0$$

which implies that also $q_n(\zeta^2) \neq 0$. Therefore, $\zeta^2$ is not a primitive $n$-th root of unity, which implies that 2 and $n$ are not relatively prime, i.e., $n$ is even. Since $n$ divides $(\dim(H))^3$ by Proposition 2.2, we see that $\dim(H)$ is also even. □

We note that the converse of the above theorem also holds: If a semisimple Hopf algebra has even dimension, it has a nontrivial self-dual simple module. To see this, look at the action of the antipode on the minimal two-sided ideals that appear in the Wedderburn decomposition. A simple module is self-dual if and only if the antipode preserves the corresponding minimal two-sided ideal. If this happens only for the one-dimensional ideal that corresponds to the trivial representation, the remaining minimal two-sided ideals can be grouped into pairs of ideals of equal dimension. As the dimension of the Hopf algebra is the sum of the dimensions of the minimal two-sided ideals, this must then be an odd number.

The arguments that we have given so far also prove two facts that are of independent interest:

**Corollary** Suppose that $H$ is a semisimple Hopf algebra over an algebraically closed field $K$ of characteristic zero.

1. If $\chi$ is an irreducible character of $H$ and $\eta$ is an irreducible character of $H^*$, then $\eta(\chi)$ is contained in the $n$-th cyclotomic field $\mathbb{Q}(\zeta_n) \subset K$, where $n$ is the exponent of $H$ and $\zeta_n$ is a primitive $n$-th root of unity of $K$.

2. If the dimension of $H$ is even, then the exponent of $H$ is also even.

**Proof.** The first statement follows from the considerations at the beginning of the proof of the theorem. The second statement hold since, if the dimension of $H$ is even, we have just seen that $H$ has a nontrivial self-dual simple module, and we have seen in the proof of the theorem that this implies that the exponent of $H$ is even. □

The second statement can be seen as a first partial answer to the question whether the exponent and the dimension of $H$ have the same prime divisors (cf. [5], Qu. 5.1, p. 138).

**5**  An important open problem in the theory of semisimple Hopf algebras is to prove that, over an algebraically closed field of characteristic zero, the dimension of a simple module divides the dimension of the Hopf algebra. This was the sixth out of a list of ten problems posed by I. Kaplansky in 1975 (cf. [9], [21]). The above theorem can be used to give a partial answer to this conjecture:

**Corollary**  Suppose that $H$ is a semisimple Hopf algebra over an algebraically closed field of characteristic zero. If $H$ has a simple module of even dimension, then the dimension of $H$ is even.

**Proof.**  Assume on the contrary that the dimension of $H$ is odd. Suppose that $\chi_1, \ldots, \chi_k$ are the irreducible characters of $H$, where $\chi_1 = \varepsilon$ is the character of the trivial module. As above, we denote the dual of an irreducible character $\chi$ by $\bar{\chi}$. From Theorem 4, we know that the trivial module is the only self-dual simple module. Therefore $k = 2l + 1$ must be odd, and we can number the characters in such a way that no pair of dual characters is contained in $\chi_2, \ldots, \chi_{l+1}$ and that the remaining characters $\chi_{l+2}, \ldots, \chi_k$ are the duals of the first:

$$\chi_{i+l} = \bar{\chi}_i$$

for $i = 2, \ldots, l+1$.

Now assume that $\chi$ is an irreducible character of even degree $n$. By Schur's lemma, the trivial character appears exactly once in the decomposition of $\chi\bar{\chi}$. This decomposition therefore has the form

$$\chi\bar{\chi} = \chi_1 + \sum_{i=2}^{l+1} m_i \chi_i + \sum_{i=l+2}^{k} m_i \chi_i$$

where $m_i \in \mathbb{N}_0$ is the multiplicity of $\chi_i$ in $\chi\bar{\chi}$. Since $\chi\bar{\chi}$ is self-dual, we must have $m_{i+l} = m_i$ for $i = 2, \ldots, l+1$. Denoting the degree of $\chi_i$ by $n_i$, we can take degrees in the above equation to get

$$n^2 = 1 + 2 \sum_{i=2}^{l+1} m_i n_i$$

Since the left hand side is even and the right hand side is odd, this is a contradiction. □

This corollary generalizes, over algebraically closed fields of characteristic zero, a result of W. D. Nichols and M. B. Richmond, who proved it in the case where the simple module has dimension 2 (cf. [19], Cor. 12, p. 306). We note that the proof also shows that the dimension of $H$ must be even if the dimension of the character ring is even.

**6** By the lifting theorems of P. Etingof and S. Gelaki, results of the type of the above theorem can be carried over to fields of positive characteristic under the additional assumption that the Hopf algebra is also cosemisimple. We therefore have the following consequence:

**Corollary** Suppose that $H$ is a semisimple cosemisimple Hopf algebra over an algebraically closed field $K$. If $H$ has a nontrivial self-dual simple module, then the dimension of $H$ is even.

**Proof.** Let us explain in detail how the lifting theorems can be applied to our situation: Obviously, we can assume that the characteristic $p$ of $K$ is positive. Since $K$ is perfect, there exists a complete discrete valuation ring $R$ of characteristic zero with residue field $K$ whose maximal ideal is generated by $p$ (cf. [20], § II.5, Thm. 3, p. 36), namely the ring of Witt vectors of $K$. Such a discrete valuation ring is unique up to isomorphism, and we denote its quotient field by $F$. By [6], Thm. 2.1, p. 855, there is an $R$-Hopf algebra $A$ with the properties that $A \otimes_R F$ is a semisimple and cosemisimple Hopf algebra over the field $F$ of characteristic zero and that $A/pA$ is isomorphic to $H$. If

$$H \cong \bigoplus_{i=1}^{k} M(n_i \times n_i, K)$$

is the Wedderburn decomposition of $H$, we can construct an isomorphism between $A$ and $\bigoplus_{i=1}^{k} M(n_i \times n_i, R)$ in such a way that the diagram

$$\begin{array}{ccc} A & \longrightarrow & \bigoplus_{i=1}^{k} M(n_i \times n_i, R) \\ \downarrow & & \downarrow \\ H & \longrightarrow & \bigoplus_{i=1}^{k} M(n_i \times n_i, K) \end{array}$$

commutes, where the vertical arrows arise from the quotient mappings from $A$ to $H$ resp. from $R$ to $K$ (cf. [13], § III.5, Lem. (5.1.16), p. 142; [14], § 22, Thm. (22.11), p. 342).

Suppose now that $V$ is a nontrivial self-dual simple $H$-module. Suppose that the index $j$ marks the minimal two-sided ideal in the above Wedderburn decomposition that corresponds to $V$, and let $e_j$ be the corresponding centrally primitive idempotent. If we assume that the trivial representation corresponds to the first indexed block, the nontriviality assertion on $V$ implies that $j \neq 1$. On the other hand, the self-duality assumption implies that the antipode preserves the minimal two-sided ideal, and therefore also the centrally primitive idempotent $e_j$. As the antipode of $A$ lifts the antipode of $H$, and there is a one-to-one correspondence between the centrally primitive idempotents of $A$ and the centrally primitive idempotents of $H$ (cf. [14], loc. cit.), the antipode of $A$ also preserves the $j$-th centrally primitive idempotent $e'_j$ of $A$. Therefore, the antipode of $A \otimes_R F$ preserves the centrally primitive idempotent $e'_j \otimes 1$, which implies that $A \otimes_R F$ has a nontrivial self-dual simple module. By Theorem 4, we see that $\dim_K H = \dim_F A \otimes_R F$ is even. □

By a similar argument, one can show the following:

**Corollary** Suppose that $H$ is a semisimple cosemisimple Hopf algebra over an algebraically closed field. If $H$ has a simple module of even dimension, then the dimension of $H$ is even.

Of course, it is also possible to use the proof of Corollary 5 directly.

Typeset using $\mathcal{AMS}$ - LaTeX